%-----------------------------------------------------------------------------
% Beginning of article.tex for opus 11
% Written in AMSTEX (not Latex).
% Inputs expm.sty which inputs AMS-J.STY for AMS journals, 
% which inputs AMSPPT.STY.
%-----------------------------------------------------------------------------
%
% AMS-TeX 2.1 sample file for journals.
%
\input amstex
\documentstyle{amsppt}
\NoBlackBoxes
\TagsOnRight
\loadbold

% Some macros
\def\quats{$\Bbb{H}$}
\def\qi{\bold i}
\def\qj{\bold j}
\def\qk{\bold k}
\def\ztwothreesix{\Bbb Z \left[ 1,\, \sqrt 2\,\qi,\, \sqrt 3\,\qj,\, \sqrt 6\, \qk \right]}

\def\localcite#1#2#3{[#2]} % author number year
% end of some macros

\settabs  12 \columns

\topmatter
\title A Quaternionic Proof of the Universality of Some Quadratic Forms\endtitle
\author Jesse Ira Deutsch \endauthor

\address Mathematics Department, University of Botswana, Private Bag 0022,
Gaborone, Botswana \endaddress

%\curraddr Department of Mathematics \endcurraddr 

\email deutschj\@mopipi.ub.bw\endemail

%  The following items provide publication information for the AMS-J logo
%\issueinfo{00}% volume number
%{0}% issue number
%{Monat}% month
%{2003}% year

%  Math Subject Classifications 
%  11P05 Waring's problem and variants, 11D57 mult. and norm form eqtns
%  11D09 quadratic and bilinear diophantine eqtns
%  11Y99 other computational number theory

\abstract 
The problem of finding all quadratic forms over $\Bbb Z$
that represent each positive integer received
significant attention in a paper of Ramanujan in 1917.
Exactly fifty four quaternary quadratic forms of this type
without cross product terms 
were shown to represent all positive integers.
The classical case of the quadratic form that is just the
sum of four squares received an alternate proof by
Hurwitz using a special ring of quaternions.
Here we prove that seven other quaternary quadratic forms
can be shown to represent all positive integers by
investigation of the corresponding quaternion rings.
\endabstract

%  \thanks will become a 1st page footnote.
%  Use \endgraf to indicate a new paragraph; a blank line or \par will
%  be recognized as an error.
%\thanks The first author was supported in part by NSF
%Grant \#000000.\endthanks

%\date Monat 00, 2002 \enddate

%\dedicatory This paper is dedicated to our advisors.\enddedicatory


\endtopmatter

\document

%------------------------------------------------
\head 1. Introduction\endhead

A generalization of the classical Four Squares Theorem
concerns the representation of all positive integers
by quadratic forms over $\Bbb Z$.
Such a quadratic form is called universal.
In a 1917 paper, Ramanujan investigated the case
of quaternary quadratic forms without cross product
terms.
It was eventually determined that there were a total
of 54 such universal quaternary quadratic forms.
Various techniques have been used to demonstrate the
universality of these forms including the theory of
ternary quadratic forms and the theory of theta functions.
See Ramanujan [10] and Fine [2] for details of these
demonstrations. 
Duke [1] is useful for background material and an alternate
approach to the problem.

In 1919 Hurwitz demonstrated the classical Four Squares
Theorem by means of a special ring of quaternions.
The Hurwitz quaternions are the $\Bbb Z$ module with
generators
$\left\{ 1,\, \qi,\, \qj ,\, {1 \over 2} (1+\qi+\qj+\qk) \right\}$.
Note that each of these generators is a unit.
It is easy to see that this set is closed under
multiplication and thus forms a non-commutative subring
of the ring of all quaternions with real coefficients.
Some of the key properties needed in the proof are
that the Hurwitz quaternions comprise a norm Euclidean ring
and that the quadratic form under consideration,
$x^2 + y^2 + z^2 + w^2$, obeys a product law.
This product law is simply an alternate statement
of the fact that the norm of a product of 
two quaternions is the product of the norms of each
quaternion.
Another useful fact is that each Hurwitz quaternion 
has an associate with rational integer coefficients
with respect to the generators
$\left\{ 1,\, \qi,\, \qj,\, \qk\, \right\}$,
though this can be avoided by using a trick of Euler.
For details see Hardy and Wright [3, chap.~XX] and Herstein [4].

It is natural to consider what other quaternary
quadratic forms on the above list of 54 can be proved
universal by consideration of appropriate subrings
of quaternions.
We first focus our consideration on those quadratic forms
which obey a product law.
One way to generate a product law is to associate a
quadratic form with an underlying $\Bbb Z$ module of quaternions
closed under multiplication
such that the norm of an element of the $\Bbb Z$ module is
the quadratic form in question.
For diagonal forms this property is
satisfied when the product of any two non-initial coefficients 
is always 
an integer square times the third. 
See Lemma 1 of section 2 for details.
In addition we create a second ring, related to the 
first, which is norm Euclidean.
This second ring is actually an order in its corresponding
quaternion algebra over $\Bbb Q$ (see Vign\'eras [13]).
Finally we show that the product of an arbitrary
element of the second ring with properly chosen units 
is an element of the $\Bbb Z$ module or a suitable submodule.

The coefficient condition is satisfied by a total of
seven forms on Ramanujan's list.
Two additional forms each have coefficients off by square integer
factors from one of the seven quadratic forms, and are amenable to
quaternionic techniques.
Using the notation $(1,\, a,\, b,\, c)$ for the
form $x^2 + a y^2 + b z^2 + c w^2$ these nine forms are
$$
\eqalign {
  &(1,\, 1,\, 1,\, 1) \quad (1,\, 1,\, 1,\, 4)
    \quad (1,\, 1,\, 2,\, 2) \quad (1,\, 1,\, 2,\, 8)
  \cr
  &(1,\, 1,\, 3,\, 3) \quad
    (1,\, 2,\, 2,\, 4) \quad (1,\, 2,\, 3,\, 6) 
    \quad (1,\, 2,\, 4,\, 8) \quad (1,\, 2,\, 5,\, 10)
    \,.
  \cr
}
\tag 1.1
$$
The first form, $(1,\, 1,\, 1,\, 1)$ corresponds to the 
classical Four Squares Theorem.
The second form $(1,\, 1,\, 1,\, 4)$ can be shown universal
using the same Hurwitz quaternions and a submodule of the
module of integer coefficient quaternions.
The forms $(1,\, 1,\, 2,\, 2)$, $(1,\, 1,\, 2,\, 8)$,
$(1,\, 2,\, 2,\, 4)$ and
$(1,\, 2,\, 4,\, 8)$ all utilize another norm Euclidean
order of quaternions in their demonstration of universality.

The forms $(1,\, 1,\, 3,\, 3)$ and $(1,\, 2,\, 3,\, 6)$
each require the introduction of new norm Euclidean orders.
Universality follows from an application of computer algebra
utilizing the appropriate units.
The last form $(1,\, 2,\, 5,\, 10)$ defies analysis by the
above technique.
The reasons why a straightforward approach will not work
in this case are quickly covered.

To fix notation, we use $\hbox{\quats}$ to represent
the ring of all real quaternions.
$$
\hbox{\quats} \ = \
  \left\{
    a_1 + a_2 \qi + a_3 \qj + a_4 \qk \,\,|\,\, 
    a_1, a_2, a_3, a_4 \,\in {\Bbb R} \right\}
\tag 1.2
$$
where $\qi,\, \qj,\, \qk$ are the standard noncommutative basis elements
with squares equal to $-1$.
We use bold type to represent quaternions, so that
the typical element of \quats\ as listed in (1.2) 
can be simply written as $\bold a$. 
The quaternion conjugate to $\bold a$ is
$$
\bold {\overline a} \ = \
    a_1 - a_2 \qi - a_3 \qj - a_4 \qk  
\tag 1.3
$$
and it is well known that $\bold a \cdot \bold {\overline a}$ is real
with zero coefficient for $\qi, \qj$ and $\qk$.
There is a canonical norm for \quats,
$$
N( \bold a ) \ = \ 
  a_1^2 + a_2^2 + a_3^2 + a_4^2
  \ = \ \bold a \cdot \bold {\overline a}
  \,. 
\tag 1.4
$$
We use the standard notation 
$\textstyle\left( {{m,\, n} \over {\Bbb Q}} \right)$
for the generalized quaternion algebra over the rationals.
With algebra generators denoted $\widehat{\bold i}$ and
$\widehat{\bold j}$ we have 
${\widehat{\bold i}\,}^2 = m$ and
${\widehat{\bold j}\,}^2 = n$ for $m,\, n \in {\Bbb Q}$ and restrict
to $m,\, n < 0$.
Thus we have a natural embedding
$\textstyle\left( {{m,\, n} \over {\Bbb Q}} \right)
  \hookrightarrow  \hbox{\quats}$
via
$$
a_1 + a_2 \, \widehat{\bold i} + a_3 \, \widehat{\bold j} 
  + a_4 \, \widehat{\bold k}
    \ \ \longrightarrow  \ \ 
  a_1 + a_2 \,\sqrt{-m} \, \qi + a_3 \, \sqrt{-n} \, \qj 
      + a_4 \,\sqrt{mn} \, \qk
\,.
\tag 1.5
$$
See Pierce [9], Scharlau [12] and Vign\'eras [13] for details.

The notation $\Bbb Z \left[ \bold a,\, \bold b,\, \dots \right]$
denotes the $\Bbb Z$ module generated by the elements
in square brackets.
$H_{n, m, \dots}$ represents a special subring or order in
the quaternions \quats.

%------------------------------------------------
\medskip
\head 2. Norm quaternionic forms\endhead

We define a quadratic form of four variables to be
norm quaternionic if it is the norm of a $\Bbb Z$ module
of quaternions that is closed under quaternion
multiplication.
More specifically we let $f(x,\, y,\, z,\, w)$ be the
quadratic form, and suppose that there is a $\Bbb Z$ 
module of quaternions with generators
$\bold g_1$, $\bold g_2$, $\bold g_3$, and $\bold g_4$
that is closed under multiplication.
Then we say that $f$ is norm quaternionic if
$$
f (x,\, y,\, z,\, w) 
    \ = \
  N (x \, \bold g_1 + y \, \bold g_2 + z \, \bold g_3 + w \, \bold g_4 )
\tag 2.1
$$
The typical example is
$x^2 + y^2 + z^2 + w^2$ which is the norm of the element
$x + y\qi + z\qj + w\qk$ of the module
$\Bbb Z \left[1, \qi, \qj, \qk \right]$.
More generally it is easy to see that
$x^2 + ay^2 + bz^2 + ab w^2$ is the norm of
$x + y \sqrt a \, \qi + z \sqrt b \, \qj + w \sqrt {ab} \, \qk$
while the elements of the module
$\Bbb Z [ 1, \sqrt a \, \qi, \sqrt b \, \qj, \sqrt {ab} \,\qk ]$
satisfy a multiplication law, as will be seen below.

Suppose we have a norm quaternionic form $f$
that represents each of two integers.
Then taking the product of the underlying quaternions
we conclude that the product of the two integers is also
represented by $f$.
Specifically we have the following.

\proclaim{Lemma 1}
Suppose $a, b, c$ are positive integers such that
$$
\sqrt {a b} \ = \ k_c \sqrt c ,
    \quad
  \sqrt {a c} \ = \ k_b \sqrt b ,
    \quad
  \sqrt {b c} \ = \ k_a \sqrt a
\tag 2.2
$$
where $k_a$, $k_b$, and $k_c$ are positive integers.
Then the quadratic form
$x^2 + a y^2 + b z^2 + c w^2$ is norm quaternionic.
\endproclaim

\demo{Proof}
We wish to show that the underlying module
$\Bbb Z [ 1, \sqrt a \,\qi, \sqrt b \,\qj, \sqrt c \,\qk ]$
is closed under quaternionic multiplication, but this is a simple exercise
in algebra since
$$
\eqalign { 
  &\left( x_1 + y_1 \sqrt a \,\qi + z_1 \sqrt b \,\qj 
         + w_1 \sqrt {c} \,\qk \right)
  \cdot
  \left( x_2 + y_2 \sqrt a \,\qi + z_2 \sqrt b \,\qj 
         + w_2 \sqrt {c} \,\qk \right)
  \cr
  &\qquad=\ 
    \left( x_1 x_2 - a \, y_1 y_2 - b \,z_1 z_2 - c \,w_1 w_2 \right)
  \cr
  &\qquad\qquad +\ \sqrt a \,\qi
    \left( x_1 y_2 + x_2 y_1 - k_a \, w_1 z_2 + k_a \, w_2 z_1 \right)
  \cr
  &\qquad\qquad + \ \sqrt b \,\qj
    \left( x_1 z_2 + x_2 z_1 + k_b \, w_1 y_2 - k_b \, w_2 y_1 \right)
  \cr
  &\qquad\qquad + \ \sqrt {c} \,\qk
    \left( x_1 w_2 + x_2 w_1 + k_c \, y_1 z_2 - k_c \, y_2 z_1 \right)
  \,.
  \qquad\hbox{\qed}
}
\tag 2.3
$$
\enddemo

Of the forms listed in $(1.1)$ all are norm quaternionic except
$(1,\, 1,\, 1,\, 4)$ and $(1,\, 1,\, 2,\, 8)$.

%------------------------------------------------
\medskip
\head 3. Orders, Norm Euclidean property \endhead

For our purposes, an order of the rational quaternion algebra
$\textstyle\left( {{m,\, n} \over {\Bbb Q}} \right)$,
or simply an order,
is a finitely generated ${\Bbb Z}$ module $H$ contained in the
algebra such that $H$ is a ring with unity and
$$
{\Bbb Q} \ \otimes_{\Bbb Z} \ H \ \simeq \ 
  \displaystyle\left( {{m,\, n} \over {\Bbb Q}} \right)
.
\tag 3.1
$$
This last condition means that the whole quaternion algebra
is spanned when the scalars on the module $H$ are extended from
$\Bbb Z$ to ${\Bbb Q}\,$.
For the general definition and properties of orders of quaternion algebras 
see Vign\'eras [13, pp.~19-21].

It is easy to observe that the intersection of an order of a rational
quaternion algebra and the real numbers is ${\Bbb Z}$.
This intersection is a finitely generated $\Bbb Z$ module,
a ring with unity, and is contained in the rationals.
If the intersection contained any rational number outside of $\Bbb Z$
it would also contain rationals with arbitrarily large denominators
when written in lowest terms.
This contradicts being a finitely generated $\Bbb Z$ module.

It is also the case that the norm map on an order of a rational 
quaternion algebra takes values in $\Bbb Z$.
This is implicit in the equivalent definitions for an order given 
in Vign\'eras [13, p.~20].
It is also not difficult to prove directly.

\proclaim {Lemma 2}
The norm of any element of an order of a rational quaternion 
algebra is a rational integer.
\endproclaim

\demo {Proof}
Let $H$ be a $\Bbb Z$ order of a rational quaternion
algebra, and let $\bold a \in H$.
By Reiner [11, p.~3] $\bold a$ is integral over 
$\Bbb Z$.
It follows from Reiner [11, p.~6]
that the minimal polynomial over $\Bbb Q$ of $\bold a$ is in
$\Bbb Z [x]$.

Note that $\bold a$ satisfies a quadratic polynomial over
$\Bbb Q$, namely,
$$
x^2 - ( \bold a + \bold {\overline a} ) \, x 
  + \bold a \, \bold {\overline a}
     \ = \ 0 \,.
\tag 3.2
$$
If this polynomial factors, since $H$ has no zero divisors
it follows that $\bold a$ must be a rational number.
By a comment just before the Lemma, it follows that 
$\bold a \in {\Bbb Z}$
so that the norm of $\bold a$ is also in $\Bbb Z$.

If the polynomial above does not factor, it is a monic polynomial
of smallest degree satisfied by $\bold a$ over $\Bbb Q$, hence it
is the minimal polynomial for $\bold a$ and the coefficients must
be in $\Bbb Z$.
\qed
\enddemo

We say that a ring $R$ of quaternions is
norm Euclidean if there is a $\delta > 0$ such that
$\forall \bold q \in \hbox{\quats}$ 
$\,\,\exists \, \bold a \in R$ with
$N (\bold q - \bold a) \leq \delta < 1$.
An associate of $r \in R$ 
is an element of the form $\varepsilon r$ or $r \varepsilon$
where $\varepsilon$ is a unit.
The next results follow using the same proofs as in the case
of Hurwitz quaternions or in the classical case of the rational
integers.
Thus we shall only mention the results, and refer to
Hardy and Wright [3] for details.

\proclaim{Lemma 3}
Let $R$ be a norm Euclidean ring of quaternions and
$\bold a, \bold b \in R$.
Suppose $\bold b \neq \bold 0$ the zero quaternion.
Then there exists $\bold g, \bold h \in R$ such that
$$
\bold a \ =\ 
  \bold g \bold b \,+\, \bold h \,,
    \qquad N( \bold h ) \, < \, N( \bold b )
\tag 3.3
$$
\endproclaim

We now need the concept of greatest common divisor.
The noncommutativity of \quats\ must be taken into
account.

\definition{Definition 4}
For any two nonzero elements $\bold a$, $\bold b$ of a subring $R$ 
of \quats, $\bold d$ is a right greatest common divisor
if it divides each of $\bold a$ and $\bold b$ on the right,
and for any $\bold e \in R$ that divides these two elements
on the right, $\bold e$ divides $\bold d$ in the right.
In notation
$$
\eqalign{
&(i) \ \exists \, \bold c_1, \, \bold c_2 \in R \ \ni \ 
   \bold a = \bold c_1 \, \bold d,\
   \bold b = \bold c_2 \, \bold d
  \cr
&(ii) \ 
  \hbox{\rm if}\ \bold e \in R \ \ni \
    \exists \, \bold f_1, \, \bold f_2 \in R,\
      \bold a = \bold f_1 \, \bold e,\
      \bold b = \bold f_2 \, \bold e
    \Longrightarrow \, \exists \, \bold g \in R \ \ni \
      \bold d = \bold g \, \bold e
    \,.
  \cr
}
\tag 3.4
$$
\enddefinition

\proclaim{Lemma 5}
The right greatest common divisor is unique up to
a factor invertible in $R$.
\endproclaim

Now we show criteria under which a right greatest common
divisor must exist.

\proclaim{Lemma 6}
Let $R$ be a ring of quaternions, $R \subset \hbox{\quats}$.
Suppose that $R$ is norm Euclidean
and that the norm on \quats, $N$, maps $R$ into the nonnegative
rational integers.
Let $\bold a$ and $\bold b$ be any nonzero elements of $R$.
Then a right greatest common divisor, $\bold d$, of 
$\bold a$ and $\bold b$
exists in $R$ and can be written as a right linear 
combination of $\bold a$ and $\bold b$, {\sl i.e.}
$$
\exists \, \bold r, \, \bold s \in R \ \, \ni \,
  \bold d \ = \
    \bold r \, \bold a \,+\, \bold s \, \bold b
\tag 3.5
$$
\endproclaim

Some information is known about norm Euclidean orders.
For orders of a ring whose field of fractions is a number field, 
the norm Euclidean property implies that the order is maximal
(see Vign\'eras [13, p.~91]).
In the case of an algebra ``totalement d\'efinie'' over
$\Bbb Q$ the norm Euclidean property implies that the reduced
discriminant is $2,\, 3$ or $5$ (see Vign\'eras [13, p.~156]).
Criteria exist concerning the Euclidean property for functions
other than the norm.
For more details on this aspect of the theory, see 
R. Markanda and V. Albis-Gonzales~[8].

%------------------------------------------------
\medskip
\head 4. The form $(1,\, 1,\, 1,\, 4)$\endhead

We start with the ring of Hurwitz quaternions
$H_{1,1,1} = \Bbb Z \, [ 1,\, \qi,\, \qj,\, 
              {1 \over 2} (1 + \qi + \qj + \qk)]$
and recall from the demonstration of the classical Four
Squares Theorem
that for every positive integer $n$ there exists a 
quaternion $\bold q \in H_{1,1,1}$ such that $N (\bold q) = n$.
Also every element of $H_{1,1,1}$ has an associate in
$\Bbb Z [1,\, \qi,\, \qj,\, \qk]$.
See Hardy and Wright [3], Herstein [4] or Hurwitz [6] for details.
An analogue of these results suffices for a demonstration of
the universality of the form $(1,\, 1,\, 1,\, 4)$ despite the
fact that this form is not itself norm quaternionic.
We set $\bold h = {1 \over 2} (1 + \qi + \qj + \qk)$ for brevity.

\proclaim{Lemma 7}
Every element of $H_{1,1,1}$ has an associate in the $\Bbb Z$-module
$\Bbb Z [1,\, \qi,\, \qj,\, 2\,\qk]$.
\endproclaim

\demo{Proof}
Write a typical $\bold a \in H_{1,1,1}$ in the form
$$
\bold a \ = \
  4 \, (c_1 \cdot 1 + c_2 \qi + c_3 \qj + c_4 \bold h )
    \,+\,
  d_1 \cdot 1 + d_2 \qi + d_3 \qj + d_4 \bold h
\tag 4.1
$$
where $c_1, \dots, c_4$ are rational integers and 
$d_1, \dots , d_4 \in \{ 0,\, 1,\, 2,\, 3 \}$.
Let 
$\bold d = d_1 + d_2\,\qi + d_3\,\qj + d_4\, \bold h$
and define $\bold c$ analogously.
Then $\bold c \in H_{1,1,1}$.
Computation shows there exists a unit $\bold u \in H_{1,1,1}$
for which 
$\bold d \cdot \bold u \in \Bbb Z [1,\, \qi,\, \qj,\, 2 \, \qk]$.
See Table I for a sample of these calculations.
Since $\bold c \cdot \bold u$ is an element of the ring
$H_{1,1,1}$ clearly $4 \, \bold c \cdot \bold u$ is an 
element of $\Bbb Z [1,\, \qi,\, \qj,\, 2 \, \qk]$.
Hence 
$$
\bold a \, \bold u \ = \
  4 \, \bold c \, \bold u + \bold d \, \bold u
\tag 4.2
$$
is also an element of $\Bbb Z [1,\, \qi,\, \qj,\, 2 \qk]$.
\qed
\enddemo

%----- start of Table I -----%
\topinsert
\topcaption{Table I}
Some associates of elements of $H_{1,1,1}$ that are members of
$\Bbb Z [ 1, \, \qi, \, \qj,  2 \, \qk ] $
\endcaption
\hfil\vbox{ \openup 0.75\jot
\halign { # \hfil  & \quad\quad # \quad\quad\hfil & \hfil #  \cr
Quaternion $\bold q$ & Unit $\bold u$ & Associate $\bold q \cdot \bold u$ \cr
\cr
$\vdots$     &    $\vdots$        &  $\vdots$
\cr
%----------------------------------------------
$ 3 +  2\,\qj +  2\,\bold h$ & $ 2\,\bold h - \qj - \qi - \, 1$ & $  4\,\qk-\,\qj+3\,\qi-1 $ \cr
$ 3 +  2\,\qj +  3\,\bold h$ & $ \bold h$ & $  2\,\qk+4\,\qj+4\,\qi-1 $ \cr
$ 3 +  3\,\qj$ & $ \, 1$ & $  3\,\qj+3 $ \cr
$ 3 +  3\,\qj +   \bold h$ & $ \bold h - \, 1$ & $  3\,\qi-4 $ \cr
$ 3 +  3\,\qj +  2\,\bold h$ & $ \qi$ & $  -4\,\qk+\,\qj+4\,\qi-1 $ \cr
$ 3 +  3\,\qj +  3\,\bold h$ & $ \bold h - \, 1$ & $  3\,\qi-6 $ \cr
$ 3 +   \qi$ & $ \, 1$ & $  \,\qi+3 $ \cr
$ 3 +   \qi +   \bold h$ & $ \bold h - \, 1$ & $  2\,\qk+\,\qj+\,\qi-3 $ \cr
$ 3 +   \qi +  2\,\bold h$ & $ \qj$ & $  2\,\qk+4\,\qj-\,\qi-1 $ \cr
$ 3 +   \qi +  3\,\bold h$ & $ \bold h - \, 1$ & $  2\,\qk+\,\qj+\,\qi-5 $ \cr
$ 3 +   \qi +   \qj$ & $ \, 1$ & $  \,\qj+\,\qi+3 $ \cr
%----------------------------------------------
$\vdots$     &    $\vdots$        &  $\vdots$
\cr
}   % end halign
} %end hbox
\hfill
%---- to the compositor -----%
%  Change these lines as you want.  It just
%  didn't look very nice without some
%  separation from the regular text.
%
\medskip
%\hrule
\bigskip
\endinsert
%----- end of Table I -----%

\proclaim{Theorem 8}
The form $(1,\, 1,\, 1,\, 4)$ is universal.
\endproclaim

\demo{Proof}
Fix any positive integer $n$.
By the proof of the classical Four Squares Theorem via
quaternions, there exists a $\bold q \in H_{1,1,1}$ 
for which $N(\bold q) = n$ (see Hardy and Wright [3]).
By the previous Lemma there exists a unit $\bold u \in H_{1,1,1}$
such that 
$\bold q \, \bold u \in \Bbb Z [1,\, \qi,\, \qj,\, 2\, \qk]$.
Write
$$
\bold q \, \bold u = x + y\, \qi + z \, \qj + w \, 2\, \qk
\tag 4.3
$$
with $x$, $y$, $z$, $w$ rational integers.
Then 
$$
n \ = \  N(\bold q) \ = \ N(\bold q \, \bold u)
  \ = \  x^2 + y^2 + z^2 + 4 \, w^2
  \,.
  \qquad\hbox{\qed}
\tag 4.4
$$
\enddemo

%------------------------------------------------
\medskip
\head 5. The forms $(1,\, 1,\, 2,\, 2)$, $(1,\, 1,\, 2,\, 8)$,
         $(1,\, 2,\, 2,\, 4)$, $(1,\, 2,\, 4,\, 8)$\endhead

We note that $(1,\, 1,\, 2,\, 8)$, $(1,\, 2,\, 2,\, 4)$ 
and $(1,\, 2,\, 4,\, 8)$
differ from the first form listed above by square factors
in the coefficients.
We call $(1,\, 1,\, 2,\, 2)$ the base form and construct
the $\Bbb Z$ module and quaternion ring for this form.
The non-base forms will then have their universality
demonstrated in a manner similar to the derivation of
the universality of $(1,\, 1,\, 1,\, 4)$ from the Hurwitz
quaternions.

Since clearly
$$
N \left( x + y \,\qi + z \sqrt 2 \,\qj + w \sqrt 2 \,\qk \right)
  \ = \ 
    x^2 + y^2 + 2 z^2 + 2 w^2
  \,.
\tag 5.1
$$
the $\Bbb Z$ module associated with this form is
$\Bbb Z \left[ 1,\, \qi,\, \sqrt 2 \,\qj,\, \sqrt 2 \, \qk \right]$.
We look for units with low denominators analogous to
Hurwitz's special quaternion 
${1 \over 2} \left(1 + \qi + \qj + \qk \right)$.
Choosing a set of these which are closed under multiplication
we have
$$
\bold v_1 \ = \ 1,\ 
  \bold v_2 \ = \ \qi,\ 
  \bold v_3 \ = \ {1 \over 2} \left( 1 + \qi + \sqrt 2 \, \qj \right),\ 
  \bold v_4 \ = \ {1 \over 2} \left( 1 + \qi + \sqrt 2 \, \qk \right)
  \,.
\tag 5.2
$$
We set 
$H_{1,2,2} = \Bbb Z \, [ \bold v_1, \bold v_2, \bold v_3, \bold v_4 ]$.
We note that 
$H_{1,2,2} \cap \Bbb R = \Bbb Z$
and the canonical norm maps $H_{1,2,2}$ into the nonnegative integers.
It is easy to see that this module is closed under quaternion
conjugation.

%----- start of Table II -----%
\topinsert
\topcaption{Table II}
Multiplication Table for $H_{1,2,2}$
\endcaption
\hfil\vbox{\offinterlineskip
\halign {\strut $#$ \ \  \vrule & \quad $#$ \quad 
& \quad\hfil $#$\hfil\quad & \quad\hfil$#$\hfil\quad 
& \quad\hfil $#$ \hfil\quad \cr
{} & \bold v_1 & \bold v_2 & \bold v_3 & \bold v_4 \cr
\noalign{\hrule}
\cr
\bold v_1 & \bold v_1 & \bold v_2 & \bold v_3 & \bold v_4  \cr
\cr
\bold v_2 & \bold v_2 & -\bold v_1 &\bold v_4 - \bold v_1 
                                       & -\bold v_3 + \bold v_2 \cr
!\cr
\bold v_3 & \bold v_3 & -\bold v_4 + \bold v_2 & \bold v_3 - \bold v_1
                                       & \bold v_2 \cr
\cr
\bold v_4 & \bold v_4 & \bold v_3 - \bold v_1
          &\bold v_4 + \bold v_3 - \bold v_2 - \bold v_1
          & \bold v_4 - \bold v_1 \cr
}   % end halign
} %end hbox
\hfill
%---- to the compositor -----%
%  Change these lines as you want.  It just
%  didn't look very nice without some
%  separation from the regular text.
%
\medskip
%\hrule
\bigskip
\endinsert
%----- end of Table II -----%

\proclaim{Lemma 9}
The $\Bbb Z$ module $H_{1,2,2}$ is closed under multiplication
and is thus an order in the quaternion algebra
$\textstyle\left( {{-1,\, -2} \over {\Bbb Q}} \right)$.
\endproclaim

\demo{Proof}
Table II, generated by computer algebra,
shows that $H_{1,2,2}$ is closed under multiplication.
It is easy to see that rational linear combinations of the
generators span all of the associated quaternion algebra.
\qed
\enddemo

\proclaim{Lemma 10}
For any rational prime $p \neq 2$ there exist rational integers
$a$ and $b$ such that 
$2 a^2 + 2 b^2 + 1 \equiv 0 \pmod{p}$.
\endproclaim

\demo{Proof}
Since $ p \neq 2$
the sets $\{ 2 a^2 \}$ and $\{ -1 - 2b^2 \}$
each contain $1 + {{p-1} \over 2}$ different elements
as $a$ and $b$ run through all the integers modulo $p$.
This is because exactly half the integers from $1$
to $p-1$ are squares modulo $p$.
If these two sets had no intersection modulo $p$ we 
would have a total of $2 \cdot {{p+1} \over 2} = p + 1$ 
different elements modulo $p$, a contradiction.
Hence there exists $a, b \in \Bbb Z$ such that
$2 a^2 \equiv -1 -2 b^2 \pmod{p}$.
\qed
\enddemo

An important step in the proof of the Four Squares Theorem
by Hurwitz quaternions is the result that any such
quaternion has an associate with integer coordinates
with respect to the basis $1,\, \qi,\, \qj,\, \qk$ 
(see Hardy and Wright [3, chap.~XX]).
The analogous result for $H_{1,2,2}$ can be obtained 
by first listing the units in the ring.

\proclaim{Lemma 11}
The units of $H_{1,2,2}$ are
$$
\pm \left\{
  1,\,\,\, \qi,\,\,\,
  {1 \over 2} \pm {1 \over 2} \qi \pm {{\sqrt 2} \over 2} \qj,\,\,\,
  {1 \over 2} \pm {1 \over 2} \qi \pm {{\sqrt 2} \over 2} \qk,\,\,\,
  {{\sqrt 2} \over 2} \qj \pm {{\sqrt 2} \over 2} \qk
\right\}
\tag 5.3
$$
or in terms of basis elements
$$
\pm \left\{
\eqalign {
  &\bold v_1,\, \bold v_2,\, 
  \bold v_3,\, \bold v_3 - \bold v_1,\, \bold v_3 - \bold v_2,\,
    \bold v_3 - \bold v_2 - \bold v_1,\,
  \bold v_4,\, \bold v_4 - \bold v_1,\,
  \cr
    &\bold v_4 - \bold v_2,\,
    \bold v_4 - \bold v_2 - \bold v_1,\,
  \bold v_4 - \bold v_3,\,
  \bold v_4 + \bold v_3 - \bold v_2 - \bold v_1
  \cr
}
\right\} \,\,.
\tag 5.4
$$
\endproclaim

\demo{Proof}
To find all the units of $H_{1,2,2}$ write a typical
element of the ring in the form
$a_1 + a_2 \qi + a_3 \sqrt 2 \, \qj + a_4 \sqrt 2 \, \qk$.
The coefficients $a_1 \dots, a_4$ are rational numbers with
denominators $1$ or $2$.
To be a unit, the sum $a_1^2 + a_2^2 + 2a_3^2 + 2a_4^2$ must
equal one.
Thus the absolute values of $a_3$ and $a_4$ are less than
or equal to
$\sqrt 2 /2$ while the absolute values of $a_1$ and $a_2$
are less than or equal to one.
$a_3$ and $a_4$ can only take the values $0$ or $\pm {1 \over 2}$
while $a_1$ and $a_2$ can take on the additional values of
$\pm 1$.

If $a_3$ and $a_4$ are both zero we have the units of the
Gaussian integers $\pm 1$, $\pm \qi$.
If $a_3$ is zero and $a_4$ is not zero then we have
$a_1^2 + a_2^2 = {1 \over 2}$
which forces both $a_1$ and $a_2$ to be nonzero.
The only solution with half integers is $a_1 = \pm {1 \over 2}$
and $a_2 = \pm {1 \over 2}$.
The situation where $a_4$ is zero and $a_3$ is not zero is completely
analogous to the previous case.
Finally if $a_3$ and $a_4$ are both nonzero then their squares
are ${1 \over 4}$ and this forces $a_1 = a_2 = 0$.
The units in this case correspond to $a_3 = \pm {1 \over 2}$
and $a_4 = \pm {1 \over 2}$.
\qed
\enddemo

\proclaim{Lemma 12}
Every element of $H_{1,2,2}$ has an associate 
in 
$\Bbb Z [1,\, \qi,\, \sqrt 2 \, \qj,\, \sqrt 2 \, \qk ]$.
\endproclaim

\demo{Proof}
Write 
$\bold a = a_1 \bold v_1 + a_2 \bold v_2 + a_3 \bold v_3 + a_4 \bold v_4$
in the form
$$
\bold a \ = \
  2 \, (c_1 \bold v_1 + c_2 \bold v_2 + c_3 \bold v_3 + c_4 \bold v_4 )
    \,+\,
  d_1 \bold v_1 + d_2 \bold v_2 + d_3 \bold v_3 + d_4 \bold v_4
\tag 5.5
$$
where $d_1, \dots , d_4 \in \{ 0,\, 1 \}$.
Let 
$\bold c = \sum_{t=1}^4 c_t \bold v_t$
and
$\bold d = \sum_{t=1}^4 d_t \bold v_t$.
Take the unit $\bold u$ from Table III such that 
$\bold d \bold u \in \Bbb Z [1,\, \qi,\, \sqrt 2 \, \qj,\, \sqrt 2 \, \qk]$.
Then
$$
\bold a \, \bold u \ = \
  2 \, \bold c \, \bold u \,+\, \bold d \, \bold u
    \ = \
  2 \, \bold {\widehat c} \,+\, \bold d \,\bold u
\tag 5.6
$$
where $\bold {\widehat c} \in H_{1,2,2}$ as $H_{1,2,2}$ is a ring.
Since twice any element of $H_{1,2,2}$ must be in
$\Bbb Z [1,\, \qi,\, \sqrt 2 \, \qj,\, \sqrt 2 \, \qk ]$
and $\bold d \, \bold u$ is also in this module, we find
that $\bold a \, \bold u$ is in the module, the result wanted.
\qed
\enddemo

%----- start of Table III -----%
\topinsert
\topcaption{Table III}
Associates of elements of $H_{1,2,2}$ that are members of
$\Bbb Z [ 1, \qi, \sqrt 2 \, \qj, \sqrt 2 \, \qk ] $
\endcaption
\hfil\vbox{ \openup 0.75\jot
\halign { # \hfil  & \quad\quad # \quad\quad\hfil & \hfil #  \cr
Quaternion $\bold q$ & Unit $\bold u$ & Associate $\bold q \cdot \bold u$ \cr
\cr
$\bold v_1$  &    $\bold v_1$     &  $1$
\cr
$\bold v_2$  &  $\bold v_2$       &  $-1 $
\cr
$\bold v_3$  &  $\bold v_4$       &  $\qi$
\cr
$\bold v_4$  &  $\bold v_3 - \bold v_2 - \bold v_1$       
                                  &  $-\qi$
\cr
$\bold v_1 + \bold v_2$  &  $\bold v_1$       &  $1 + \qi$
\cr
$\bold v_1 + \bold v_3$  &  $\bold v_3$       &  $\qi + \sqrt 2 \, \qj$
\cr
$\bold v_1 + \bold v_4$  &  $\bold v_4$       &  $\qi + \sqrt 2 \, \qk$
\cr
$\bold v_2 + \bold v_3$  &  $\bold v_4 - \bold v_1$       
                                              &  $-1 - \sqrt 2 \, \qj$
\cr
$\bold v_2 + \bold v_4$  &  $\bold v_3 - \bold v_2$       
                                              &  $1 + \sqrt 2 \, \qk$
\cr
$\bold v_3 + \bold v_4$  &  $\bold v_4 - \bold v_3$       
                                              &  $\qi - \sqrt 2 \, \qj$
\cr
$\bold v_1 + \bold v_2 + \bold v_3$  &  $\bold v_4 - \bold v_2 - \bold v_1$       
                         &  $-\qi - \sqrt 2 \, \qj + \sqrt 2 \, \qk$
\cr
$\bold v_1 + \bold v_2 + \bold v_4$  &  $\bold v_3$       
                         &  $\qi + \sqrt 2 \, \qj + \sqrt 2 \, \qk$
\cr
$\bold v_1 + \bold v_3 + \bold v_4$  &  $\bold v_4$       
                         &  $2\qi + \sqrt 2 \, \qk$
\cr
$\bold v_2 + \bold v_3 + \bold v_4$  &  $\bold v_3$       
                         &  $-1 + \qi + \sqrt 2 \, \qj + \sqrt 2 \, \qk$
\cr
$\bold v_1 + \bold v_2 + \bold v_3 + \bold v_4$  
                         &  $\bold v_4 + \bold v_3 - \bold v_2 - \bold v_1$       
                         &  $-1 + 2 \sqrt 2 \, \qk$
\cr
}   % end halign
} %end hbox
\hfill
%---- to the compositor -----%
%  Change these lines as you want.  It just
%  didn't look very nice without some
%  separation from the regular text.
%
\medskip
%\hrule
\bigskip
\endinsert
%----- end of Table III -----%

Any number of the form $\bold u_1 \bold a \bold u_2$
with $\bold u_1$ and $\bold u_2$ units
will be called a two sided associate of $\bold a$.

\proclaim{Lemma 13}
Every element of $H_{1,2,2}$ has a two sided associate 
in each of the modules
$$
\Bbb Z [1,\, \qi,\, \sqrt 2 \, \qj,\, 2\, \sqrt 2 \, \qk ],
  \quad
\Bbb Z [1,\, 2\, \qi,\, \sqrt 2 \, \qj,\, \sqrt 2 \, \qk ],
  \quad
\Bbb Z [1,\, 2\, \qi,\, \sqrt 2 \, \qj,\, 2 \,\sqrt 2 \, \qk ]
\,.
\tag 5.7
$$
\endproclaim

\demo{Proof}
The proof is similar to the previous demonstration.
Write 
$\bold a = a_1 \bold v_1 + a_2 \bold v_2 + a_3 \bold v_3 + a_4 \bold v_4$
in the form
$$
\bold a \ = \
  4 \, (c_1 \bold v_1 + c_2 \bold v_2 + c_3 \bold v_3 + c_4 \bold v_4 )
    \,+\,
  d_1 \bold v_1 + d_2 \bold v_2 + d_3 \bold v_3 + d_4 \bold v_4
\tag 5.8
$$
where $d_1, \dots , d_4 \in \{ 0,\, 1,\, 2,\, 3 \}$.
Let 
$\bold c = \sum_{t=1}^4 c_t \bold v_t$
and
$\bold d = \sum_{t=1}^4 d_t \bold v_t$.
Computation shows that there exists units
$\bold u_1$ and $\bold u_2$ in $H_{1,2,2}$ such that
$\bold u_1 \bold d \bold u_2 \in 
   \Bbb Z [1,\, \qi,\, \sqrt 2 \, \qj,\, 2 \,\sqrt 2 \, \qk]$.
We have 
$$
\bold u_1 \, \bold a \, \bold u_2 \ = \
  4 \, \bold u_1 \, \bold c \, \bold u_2 
      \,+\, \bold u_1 \bold d \, \bold u_2
    \ = \
  4 \, \bold {\widehat c} \,+\, \bold u_1 \, \bold d \,\bold u_2
\tag 5.9
$$
where $\bold {\widehat c} \in H_{1,2,2}$ as $H_{1,2,2}$ is a ring.
Since four times any element of $H_{1,2,2}$ must be in
$\Bbb Z [2,\, 2\,\qi,\, 2\,\sqrt 2 \, \qj,\, 2\,\sqrt 2 \, \qk ]$
we find $\bold u_1 \, \bold a \, \bold u_2$ is in
$\Bbb Z [1,\, \qi,\, \sqrt 2 \, \qj,\, 2 \,\sqrt 2 \, \qk]$.

A similar argument holds in the other two cases of the Lemma.
\qed
\enddemo

\proclaim{Lemma 14}
$H_{1,2,2}$ is norm Euclidean.
\endproclaim

\demo{Proof}
Fix any $\bold q \in \hbox{\quats}$.
We may choose $a_4 \bold v_4$ such that $a_4 \in \Bbb Z$
and the absolute value of the coefficient of $\qk$ of
$\bold q - a_4 \bold v_4$ is less than or equal to 
$\sqrt 2 /4$.
Similarly we find $a_3 \bold v_3$ with $a_3 \in \Bbb Z$
such that $\bold q -a_3 \bold v_3 -a_4 \bold v_4$
additionally has the absolute value of the coefficient
of $\qj$ less than or equal to $\sqrt 2 /4$.

Then we need only subtract off a further $a_1 + a_2 \qi$
with rational integral $a_1$, $a_2$ to force the coefficients
of $\qi$ and $1$ to be less than $1/2$ in absolute value.
Consequently
$$
N ( \bold q - a_1 \bold v_1 -a_2 \bold v_2 -a_3 \bold v_3
            - a_4 \bold v_4 )
  \ \leq \
    2 \,\left( {1 \over 2} \right)^2
      \,+\, 
    2 \, \left( {{\sqrt 2} \over 4} \right)^2 
  \ = \ 
    {3 \over 4}
  \,.
\tag 5.10
$$
Thus the norm Euclidean property holds.
\qed
\enddemo

\proclaim{Theorem 15}
Every rational prime $p$ can be represented by the form
$(1,\, 1,\, 2,\, 2)$.
\endproclaim

\demo{Proof}
If $p=2$ this is trivial.
Otherwise, by Lemma 10 there exists $a,b \in \Bbb Z$ such that
$2a^2 + 2b^2 + 1 \equiv 0 \pmod {p}$.
We may choose $|a|, |b| \leq (p-1)/2$.
Thus
$$
\eqalign{
    2a^2 + 2b^2 + 1 
  \ &\leq \
    2 \, {{(p-1)^2} \over 4} +
      2 \, {{(p-1)^2} \over 4} + 1
  \cr
  \ &= \ 
    (p-1)^2 + 1 \ = \ p^2 - 2p + 2 \ < \ p^2
  \cr
}
\tag 5.11
$$
for all $p > 1$.
Hence $2a^2 + 2b^2 + 1 = p r$ where $0 < r < p$.
Set $\bold a = 1 - a \sqrt 2 \, \qj + b \sqrt 2 \, \qk$
and note $N( \bold a ) = 2a^2 + 2b^2 + 1$.
Calculate
$$
\eqalign{
2b \, \bold v_4 \,-\, &2a \, \bold v_3 \,+\, \bold v_1 
    \,-\, (b-a) \, (\bold v_1 + \bold v_2 )
  \cr
  \ &= \ 
    b \, (1 + \qi + \sqrt 2 \, \qk) - a \, (1 + \qi + \sqrt 2 \, \qj)
    + 1 -(b-a) \, (1 + \qi)
  \cr
  \ &= \
    b \sqrt 2 \, \qk - a \sqrt 2 \, \qj + 1
    \ = \ \bold a
  \,.
  \cr  
}
\tag 5.12
$$
So $\bold a \in H_{1,2,2}$.
Since $H_{1,2,2}$ is norm Euclidean, right greatest common
divisors exist.
Suppose the right greatest common divisor of $\bold a$ and $p$
is a unit $\bold u$.
Then there exists $\bold s, \bold t \in H_{1,2,2}$ such
that 
$$
\bold u \ = \ \bold s \, \bold a \,+\, \bold t \, p
  \,.
\tag 5.13
$$
Taking conjugates we find
$$
\bold {\overline u} \ = \  
     \bold {\overline a} \, \bold {\overline s} 
   \,+\, \bold {\overline t} \, p
\tag 5.14
$$
so by multiplying the previous two equations we have
for some $\bold q \in H_{1,2,2}$
$$
\eqalign{
1 \ = \ \bold u \, \bold {\overline u}
  \ &= \
      \bold s \, \bold a \, \bold {\overline a} \, \bold {\overline s} 
      \,+\, p \, \bold q
  \cr
  \ &= \
      N(\bold s) \, N(\bold a ) \,+\, p \, \bold q
  \cr
  \ &= \
      p \, r \, N(\bold s) \,+\, p \, \bold q
  \cr
}
\tag 5.15
$$
Here we use the fact that $H_{1,2,2}$ is closed under
conjugation.
Dividing through by $p$ implies
$$
{1 \over p} \ = \
    r \, N(\bold s) \,+\, \bold q 
  \ \in \ 
   H_{1,2,2}
\tag 5.16
$$
but $H_{1,2,2} \cap \Bbb R = \Bbb Z$, a contradiction.
Hence the right greatest common divisor of $\bold a$ and
$p$ is not a unit.
Call this right greatest common divisor $\bold d$.
We may write
$$
\bold a \ = \ \bold f \, \bold d,
    \quad
  p \ = \ \bold g \, \bold d,
    \qquad
  N(\bold d) \ > \ 1
\tag 5.17
$$
for some $\bold f, \, \bold g \ \in \ H_{1,2,2}$.
But then
$p^2  =  N(p)  =  N(\bold g) \, N(\bold d)$
so $N(\bold d)$ must divide $p^2$.
Also
$p r = N(\bold a) = N(\bold f) \, N(\bold d)$
so $N(\bold d)$ divides $p r$.
Thus $N(\bold d)$ divides the greatest common divisor of
$p^2$ and $pr$ in the ring $\Bbb Z$, which is simply $p$.
Since $N(\bold d)$ cannot equal one, it must equal $p$.

Choose an associate of $\bold d$ that is in the module
$\Bbb Z [1,\, \qi,\, \sqrt 2 \, \qj,\, \sqrt 2 \, \qk ]$.
Denote this associate $\bold {\widehat d}$.
Let
$$
\bold {\widehat d} \ = \
  a_1 + a_2 \qi + a_3 \sqrt 2 \, \qj + a_4 \sqrt 2 \, \qk ,
    \qquad
  a_1, a_2, a_3, a_4 \in \Bbb Z
    \,.
\tag 5.18
$$
Then 
$$
p \ = \ N(\bold d) \ = \ N(\bold {\widehat d})
    \ = \
  a_1^2 + a_2^2 + 2a_3^2 + 2 a_4^2
  \,.
\tag 5.19
$$
\enddemo

\proclaim{Theorem 16}
Every positive integer can be represented by the form
$(1,\, 1,\, 2,\, 2)$.
\endproclaim

\demo{Proof}
This follows immediately as the quadratic form in
question is norm quaternionic and also represents 
the integer one.
\qed
\enddemo

\proclaim{Theorem 17}
Every positive integer can be represented by each of
the forms
$(1,\, 1,\, 2,\, 8)$,
$(1,\, 2,\, 2,\, 4)$,
$(1,\, 2,\, 4,\, 8)$.
\endproclaim

\demo{Proof}
The proof for the case $(1,\, 1,\, 2,\, 2)$ showed that
for each positive integer $n$ there exists 
$\bold q \in H_{1,2,2}$ such
that $N(\bold q) = n$.
For the first quadratic form of the Theorem, choose
a two sided associate of $\bold q$ that is in the module
$\Bbb Z [1,\, \qi,\, \sqrt 2 \, \qj,\, 2\, \sqrt 2 \, \qk ]$.
Then there exist units $\bold u_1,\, \bold u_2$ in 
$H_{1,2,2}$ such that $\bold u_1 \, \bold q \, \bold u_2$
is in this module.
It follows that
$$
\eqalign{
n \ &= \ N( \bold q) \ = \ 
  N (\bold u_1 \, \bold q \, \bold u_2)
    \ = \
  N( x + y \qi + z \, \sqrt 2 \, \qj + w \, 2 \, \sqrt 2 \, \qk)
    \cr
    \ &= \ 
  x^2 + y^2 + 2 \, z^2 + 8 \, w^2
    \cr
}
\tag 5.20
$$
for integers $x,\, y,\, z,\, w$.
Similar argumentation demonstrates the Theorem in the case of 
the other two quadratic forms.

Alternatively, one may proceed from the representation of
rational primes $p$ to universality directly in the case of
the norm Euclidean forms
$(1,\, 2,\, 2,\, 4)$ and $(1,\, 2,\, 4,\, 8)$.
By the argument in the previous paragraph, each rational prime and $1$
is represented by each of the above two norm Euclidean forms.
Using the product law for norm Euclidean forms, universality follows
immediately.
\qed
\enddemo

\medskip
%------------------------------------------------
\head 6. The form $(1,\, 2,\, 3,\, 6)$\endhead

The proof follows closely that of the form
$(1,\, 1,\, 2,\, 2)$.
Considering the identity
$$
N \left( x + y \sqrt 2 \, \qi + z \sqrt 3 \, \qj + w \sqrt 6 \, \qk \right)
  \ = \ 
    x^2 + 2\, y^2 + 3 \, z^2 + 6\, w^2
  \,.
\tag 6.1
$$
the $\Bbb Z$ module associated with this form is
$\Bbb Z \left[ 1,\, \sqrt 2\,\qi,\, \sqrt 3\,\qj,\, \sqrt 6\, \qk \right]$.
Again we look for likely units and find
$$
\displaylines{
    \bold w_1 \ = \ 1,\ \ 
    \bold w_2 \ = \ {1 \over 2} + {{\sqrt 3} \over 2} \qj,\ \
    \bold w_3 \ = \ {1 \over 2} + {{\sqrt 3} \over 6} \qj
                    + {{\sqrt 6} \over 3} \qk,
  \cr
   \hfill (6.2)
  \cr
    \bold w_4 \ = \ {{\sqrt 2} \over 2}\qi + {{\sqrt 3} \over 3}\qj
                    + {{\sqrt 6} \over 6}\qk
    \,.
  \cr
}
$$
We set 
$H_{2,3,6} = \Bbb Z \, [ \bold w_1, \bold w_2, \bold w_3, \bold w_4 ]$.

\proclaim{Lemma 18}
The $\Bbb Z$ module $H_{2,3,6}$ is closed under multiplication
and is thus an order in the quaternion algebra
$\textstyle\left( {{-2,\, -3} \over {\Bbb Q}} \right)$.
\endproclaim

\demo{Proof}
Table IV, generated by computer algebra,
demonstrates closure under multiplication.
The proof follows as in Lemma 9 above.
\qed
\enddemo

%----- start of Table IV -----%
\topinsert
\topcaption{Table IV}
Multiplication Table for $H_{2,3,6}$
\endcaption
\hfil\vbox{\offinterlineskip
\halign {\strut $#$ \ \  \vrule & \quad $#$ \quad 
& \quad\hfil $#$\hfil\quad & \quad\hfil$#$\hfil\quad 
& \quad\hfil $#$ \hfil\quad \cr
{} & \bold w_1 & \bold w_2 & \bold w_3 & \bold w_4 \cr
\noalign{\hrule}
\cr
\bold w_1 & \bold w_1 & \bold w_2 & \bold w_3 & \bold w_4  \cr
\cr
\bold w_2 & \bold w_2 & \bold w_2 - \bold w_1 &\bold w_4  
                                   & \bold w_4 - \bold w_3  \cr
\cr
\bold w_3 & \bold w_3 & -\bold w_4 + \bold w_3 + \bold w_2 - \bold w_1 
          & \bold w_3 - \bold w_1  & \bold w_2 - \bold w_1 \cr
\cr
\bold w_4 & \bold w_4 & \bold w_3 - \bold w_1
          & \bold w_4 - \bold w_2  & - \bold w_1 \cr
}   % end halign
} %end hbox
\hfill
%---- to the compositor -----%
%  Change these lines as you want.  It just
%  didn't look very nice without some
%  separation from the regular text.
%
\medskip
%\hrule
\bigskip
\endinsert
%----- end of Table IV -----%

It follows that
$H_{2,3,6} \cap \Bbb R = \Bbb Z$
while the canonical norm maps $H_{2,3,6}$ into the nonnegative integers.
It is also easy to see that this module is closed under quaternion
conjugation.

\proclaim{Lemma 19}
For any rational prime $p$ there exist rational integers
$a$ and $b$ such that 
$a^2 + 2 b^2 + 3 \equiv 0 \pmod{p}$.
\endproclaim

\demo{Proof}
If $p = 2$ this is trivial.
Suppose $p \ne 2$.
The result follows from the same line of reasoning
as in Lemma~10 using the sets
$\{ a^2 \}$ and $\{ -3 - 2b^2 \}$ instead.
\qed
\enddemo

\proclaim{Lemma 20}
The units of $H_{2,3,6}$ are
$$
\pm \left\{
\eqalign {
  &\bold w_1,\, \bold w_2,\, \bold w_3,\, \bold w_4,\,
    \bold w_4 - \bold w_3,\, \bold w_3 - \bold w_2,\,
  \cr 
  &\bold w_2 - \bold w_1,\, \bold w_3 - \bold w_1,\,
    \bold w_4 - \bold w_2,\,
  \cr
  &\bold w_4 - \bold w_2 + \bold w_1,\,
    \bold w_4 - \bold w_3 - \bold w_2 + \bold w_1,\,
    \bold w_4 - \bold w_3 + \bold w_1
  \cr
}
\right\} \,\,.
\tag 6.3
$$
\endproclaim

\demo{Proof}
Write a typical element of $H_{2,3,6}$ in the form
$\bold w = a_1 \bold w_1 + a_2 \bold w_2 +a_3 \bold w_3 + a_4 \bold w_4$
with $a_1, \dots , a_4 \in \Bbb Z$.
The coefficient of $\qi$ in $\bold w$ is $a_4 {{\sqrt 2} \over 2}$.
To be a unit the sums of the squares of the $\qi$, $\qj$, $\qk$ and
real coefficients must equal one.
Thus ${1 \over 2} \,  a_4^2 \,\leq\, 1$.
This forces $| a_4 | \leq 1$.

For the $\qk$ coefficient of $\bold w$ we have
$a_4 {{\sqrt 6} \over 6} + a_3 {{\sqrt 6} \over 3}
  =  (a_4 + 2 a_3 ){{\sqrt 6} \over 6}$,
the square of which must be less than or equal to one.
Taking square roots, we have 
$ |a_4 + 2 a_3| \leq \sqrt 6$ so
$$
\eqalign {
    &2 | a_3 | - | a_4 | \ \leq \ |a_4 + 2 a_3 | \ \leq \ \sqrt 6
  \cr
    \Longrightarrow \quad &|a_3| \ \leq \ {1 \over 2} (\sqrt 6 + 1)
    \ = \ 1.72 \dots
}
\tag 6.4
$$
which implies that $| a_3 | \leq 1$.
Similar reasoning with the $\qj$ coefficient yields 
$|a_2| \leq 2$ while working with the real coefficient
gives $|a_1| \leq 2$.

A computer scan was run using a computer algebra system
to find all the units, which are listed in (6.3).
\qed
\enddemo

Computation implies that not every element of $H_{2,3,6}$
has an associate which is an element of the module
$\ztwothreesix$.
However, the following weaker statement suffices for our
purposes.

\proclaim{Lemma 21}
Every element of $H_{2,3,6}$ has a two sided associate in
the module
$\ztwothreesix$.
\endproclaim

\demo{Proof}
Write 
$\bold a = a_1 \bold w_1 + a_2 \bold w_2 + a_3 \bold w_3 + a_4 \bold w_4$
in the form
$$
\bold a \ = \
  6 \, (c_1 \bold w_1 + c_2 \bold w_2 + c_3 \bold w_3 + c_4 \bold w_4 )
    \,+\,
  d_1 \bold w_1 + d_2 \bold w_2 + d_3 \bold w_3 + d_4 \bold w_4
\tag 6.5
$$
where $d_1, \dots , d_4 \in \{ 0,\, 1,\, 2,\, 3,\, 4,\,  5 \}$.
Let 
$\bold c = \sum_{t=1}^4 c_t \bold w_t$
and
$\bold d = \sum_{t=1}^4 d_t \bold w_t$.
Computation shows that for every possible value of
$\bold d$ there exist units $\bold u_1$, $\bold u_2$
in $H_{2,3,6}$ such that $\bold u_1 \, \bold d \, \bold u_2$ 
is an element of $\ztwothreesix$.
Table V contains a section of the relevant computer
generated output.

Consequently
$$
\bold u_1 \, \bold a \, \bold u_2 \ = \
  6 \, \bold u_1 \, \bold c \, \bold u_2  
       + \bold u_1 \, \bold d \, \bold u_2
    \ = \
  6 \, \bold {\widehat c} + \bold u_1 \, \bold d \, \bold u_2
\tag 6.6
$$
where $\bold {\widehat c} \in H_{2,3,6}$ as $H_{2,3,6}$ is a ring.
Since six times any element of $H_{2,3,6}$ must be in
$\ztwothreesix$
and $\bold u_1 \, \bold d \, \bold u_2$ is also in this module, we find
that $\bold u_1 \, \bold a \, \bold u_2$ is in the module, 
the result wanted.
\qed
\enddemo

% Any number of the form $\bold u_1 \bold a \bold u_2$
% with $\bold u_1$ and $\bold u_2$ units
% will be called a two sided associate of $\bold a$.

%----- start of Table V -----%
\topinsert
\topcaption{Table V}
Some two sided associates of elements of $H_{2,3,6}$ that are members of
$\ztwothreesix$
\endcaption
\hfil\vbox{ \openup 0.75\jot
\halign { # \quad\hfil & # \hfil  & \quad\quad # \quad\quad\hfil & \hfil # \cr
Unit $\bold u_1$ & Quaternion $\bold q$ & Unit $\bold u_2$ & $\bold u_1 \cdot \bold q \cdot \bold u_2$ \cr
\cr
$\vdots$    & $\vdots$     &  $\vdots$          &  $\vdots$\cr
%-----------------------------------------
$ \bold w_3$ & $  2\,\bold w_2 +   \bold w_3 +  4\,\bold w_4$ & $ \bold w_4 - \bold w_2$ & $  -2\sqrt 6 \,\qk+\sqrt 2 \,\qi+3 $ \cr
$ \bold w_3 - \bold w_4$ & $  2\,\bold w_2 +   \bold w_3 +  5\,\bold w_4$ & $ \bold w_3 - \bold w_1$ & $  \sqrt 6 \,\qk+4\sqrt 2 \,\qi-3 $ \cr
$ \bold w_1$ & $  2\,\bold w_2 +  2\,\bold w_3$ & $ \bold w_3$ & $ 			   \sqrt 6 \,\qk+\sqrt 3 \,\qj+\sqrt 2 \,\qi-1 $ \cr
$ \bold w_3$ & $  2\,\bold w_2 +  2\,\bold w_3 +   \bold w_4$ & $ \bold w_2$ & $  -2\sqrt 2 \,\qi-3 $ \cr
$ \bold w_1$ & $  2\,\bold w_2 +  2\,\bold w_3 +  2\,\bold w_4$ & $ \bold w_1$ & $ 		       \sqrt 6 \,\qk+2\sqrt 3 \,\qj+\sqrt 2 \,\qi+2 $ \cr
$\vdots$   & $\vdots$                 &  $\vdots$          &  $\vdots$
\cr
%-------------------------------------------
$ \bold w_1$ & $  \bold w_1 +   \bold w_2 +  3\,\bold w_4$ & $ \bold w_3$ & $  \sqrt 6 \,\qk+2\sqrt 2 \,\qi-1 $ \cr
$ \bold w_3 - \bold w_4$ & $  \bold w_1 +   \bold w_2 +  4\,\bold w_4$ & $ \bold w_1$ & $  2\sqrt 3 \,\qj-\sqrt 2 \,\qi+3 $ \cr
$ \bold w_2$ & $  \bold w_1 +   \bold w_2 +  5\,\bold w_4$ & $ \bold w_3 - \bold w_1$ & $  -3\sqrt 3 \,\qj+\sqrt 2 \,\qi+2 $ \cr
$ \bold w_1$ & $  \bold w_1 +   \bold w_2 +   \bold w_3$ & $ \bold w_4$ & $  \sqrt 3 \,\qj+\sqrt 2 \,\qi-1 $ \cr
$ \bold w_1$ & $  \bold w_1 +   \bold w_2 +   \bold w_3 +   \bold w_4$ & $ \bold w_4$ & $  \sqrt 3 \,\qj+\sqrt 2 \,\qi-2 $ \cr
%--------------------------------------
$\vdots$   & $\vdots$                 &  $\vdots$          &  $\vdots$
\cr
}   % end halign
} %end hbox
\hfill
%---- to the compositor -----%
%  Change these lines as you want.  It just
%  didn't look very nice without some
%  separation from the regular text.
%
\medskip
%\hrule
\bigskip
\endinsert
%----- end of Table V -----%

\proclaim{Lemma 22}
$H_{2,3,6}$ is norm Euclidean.
\endproclaim

\demo{Proof}
Fix any $\bold q \in \hbox{\quats}$.
We may choose $a_4 \bold w_4$ such that $a_4 \in \Bbb Z$
and the absolute value of the coefficient of $\qi$ of
$\bold q - a_4 \bold w_4$ is less than or equal to 
$\sqrt 2 /4$.
Similarly we find $a_3 \bold w_3$ with $a_3 \in \Bbb Z$
such that $\bold q -a_3 \bold w_3 -a_4 \bold w_4$
additionally has the absolute value of the coefficient
of $\qk$ less than or equal to $\sqrt 6 /6$.

Now we find an integer multiple of $\bold w_2$
so that the $\qj$ coefficient of 
$\bold q - a_2 \bold w_2 - a_3 \bold w_3 - a_4 \bold w_4$
is additionally less than $\sqrt 3 / 4$ in absolute value.
We need only adjust the above difference by a rational
integer $a_1 \bold w_1$ to force the real part less than 
one half in absolute value.
All together we have 
$$
N \left( \bold q - \sum_{n=1}^4 \, a_n \bold w_n  \right)
  \, \leq \,
    \left( {1 \over 2} \right)^2
      + 
    \left( {{\sqrt 2} \over 4} \right)^2
      +
    \left( {{\sqrt 3} \over 4} \right)^2
      +
    \left( {{\sqrt 6} \over 6} \right)^2 
  \, = \,
    {{35} \over {48}}
  \,.
\tag 6.7
$$
Thus the norm Euclidean property holds.
\qed
\enddemo

\proclaim{Theorem 23}
Every rational prime $p$ can be represented by the form
$(1,\, 2,\, 3,\, 6)$.
\endproclaim

\demo{Proof}
For $p = 2$ this is trivial.
For odd primes $p$,
by Lemma 19 there exists $a,b \in \Bbb Z$ such that
$a^2 + 2b^2 + 3 \equiv 0 \pmod {p}$.
By reasoning similar to the previous case
we show that there exist integers $a$ and $b$
for which $a^2 + 2b^2 + 3 = pr$ where $0 < r < p$.
Setting $\bold a = a - b \sqrt 2 \, \qi + \sqrt 3 \, \qj$
it follows that $N(\bold a) = a^2 + 2b^2 + 3 = pr$.

Some algebra shows
$$
\eqalign{
  (a - b - 1) &\bold w_1 + (b +2) \bold w_2 + b \,\bold w_3 - 2b \,\bold w_4
  \cr
%  &= \
%    a - b - 1 + 
%       (b + 2) \left( {1 \over 2} + {{\sqrt 3} \over 2} \qj \right)
%      + b \left( {1 \over 2} + {{\sqrt 3} \over 6} \qj 
%          + {{\sqrt 6} \over 3} \qk \right)
%  \cr
%      &\qquad - 2b \left( {{\sqrt 2} \over 2}\qi + {{\sqrt 3} \over 3}\qj
%                    + {{\sqrt 6} \over 6}\qk \right)
%  \cr
  &= \
    a - b \sqrt 2 \, \qi + \sqrt 3 \, \qj
  \,.
}
\tag 6.8
$$
Thus $\bold a \in H_{2,3,6}$ and the proof proceeds as
in the previous case until the issue of associates arises.
As before we find $\bold d \in H_{2,3,6}$ such that
$N(\bold d) = p$.

Choose units $\bold u_1$ and $\bold u_2$
such that $\bold u_1 \, \bold d \bold u_2$ is in the module
$\ztwothreesix$.
%$\Bbb Z [1,\, \qi,\, \sqrt 2 \, \qj,\, \sqrt 2 \, \qk ]$.
%Denote this associate $\bold {\widehat d}$.
Let
$$
\bold u_1 \bold d \bold u_2 \ = \
  a_1 + a_2 \sqrt 2 \,\qi + a_3 \sqrt 3 \,\qj + a_4 \sqrt 6 \, \qk ,
    \qquad
  a_1, a_2, a_3, a_4 \in \Bbb Z
    \,.
\tag 6.9
$$
Then 
$$
p \ = \ N(\bold d) \ = \ N(\bold u_1 \, \bold d \, \bold u_2 )
    \ = \
  a_1^2 + 2a_2^2 + 3a_3^2 + 6 a_4^2
  \,.
\tag 6.10
$$
\enddemo

\proclaim{Theorem 24}
Every positive integer can be represented by the form
$(1,\, 2,\, 3,\, 6)$.
\endproclaim

\demo{Proof}
Again, this follows immediately as the quadratic form in
question is norm quaternionic and also represents the
integer one.
\qed
\enddemo

\medskip
%------------------------------------------------
\head 7. The form $(1,\, 1,\, 3,\, 3)$\endhead

The demonstration of the universal property of
the form $(1,\, 1,\, 3,\, 3)$ follows the previous cases
closely with one small difficulty.
Two sided units only reduce matters so far, and
we need an extension of a trick of Euler to transform
a representation of $2 \, n$ to one of $n$.
Set
$$
\displaylines{
    \bold x_1 \ = \ 1,\ \ 
    \bold x_2 \ = \ {1 \over 2} \qi - {{\sqrt 3} \over 2} \, \qk,\ \
    \bold x_3 \ = \ {{1} \over 4} \, \qi + {{\sqrt 3} \over 2} \, \qj
                    - {{\sqrt 3} \over 4} \, \qk,
  \cr
   \hfill (7.1)
  \cr
    \bold x_4 \ = \ {1 \over 2}  + {3 \over 4} \, \qi
                    + {{\sqrt 3} \over 4} \, \qk
  \cr
}
$$
and let $H_{1,3,3} = \Bbb Z [\bold x_1,\, \bold x_2,\, \bold x_3,\, 
                     \bold x_4]$.

\proclaim{Lemma 25}
The $\Bbb Z$ module $H_{1,3,3}$ is closed under multiplication
and is thus an order in the quaternion algebra
$\textstyle\left( {{-1,\, -3} \over {\Bbb Q}} \right)$.
\endproclaim

\demo{Proof}
Table VI, generated by computer algebra,
demonstrates closure under multiplication.
Consideration of rational linear combinations of the
basis elements completes the proof.
\qed
\enddemo

%----- start of Table VI -----%
\topinsert
\topcaption{Table VI}
Multiplication Table for $H_{1,3,3}$
\endcaption
\hfil\vbox{\offinterlineskip
\halign {\strut $#$ \ \  \vrule & \quad $#$ \quad 
& \quad\hfil $#$\hfil\quad & \quad\hfil$#$\hfil\quad 
& \quad\hfil $#$ \hfil\quad \cr
{} & \bold x_1 & \bold x_2 & \bold x_3 & \bold x_4 \cr
\noalign{\hrule}
\cr
\bold x_1 & \bold x_1 & \bold x_2  & \bold x_3 & \bold x_4  \cr
\cr
\bold x_2 & \bold x_2 & -\bold x_1 & \bold x_4 - \bold x_1
                                   & -\bold x_3 + \bold x_2  \cr
\cr
\bold x_3 & \bold x_3 & -\bold x_4 &  - \bold x_1  
                                   & \bold x_2  \cr
\cr
\bold x_4 & \bold x_4 & \bold x_3  & \bold x_3 - \bold x_2  
                                   & \bold x_4 - \bold x_1 \cr
}   % end halign
} %end hbox
\hfill
%---- to the compositor -----%
%  Change these lines as you want.  It just
%  didn't look very nice without some
%  separation from the regular text.
%
\medskip
%\hrule
\bigskip
\endinsert
%----- end of Table VI -----%

Just as for the previously considered $\Bbb Z$ modules $H_{a,b,c}$
we note that
$H_{1,3,3} \cap \Bbb R = \Bbb Z$.
Further the canonical norm on quaternions maps $H_{1,3,3}$
into the nonnegative rational integers.

\proclaim{Lemma 26}
The units in $H_{1,3,3}$ are
$$
\pm\left\{ 
  \bold x_1,\, \bold x_2,\, \bold x_3,\, \bold x_4,\, 
     \bold x_4 - \bold x_1,\, \bold x_3 - \bold x_2
  \right\}
\tag 7.2
$$
\endproclaim

\demo{Proof}
We wish to solve
$$
1 \ = \
  N (a_1 \bold x_1 + a_2 \bold x_2 + a_3 \bold x_3 + a_4 \bold x_4 )
    \ = \ 
  a_1^2 + a_2^2 + a_3^2 + a_4^2 + a_1 a_4 + a_2 a_3
\tag 7.3
$$
or equivalently
$$
\eqalign{
2 \ &= \
  2 a_1^2 + 2 a_2^2 + 2 a_3^2 + 2 a_4^2 + 2 a_1 a_4 + 2 a_2 a_3
    \cr
    \ &= \ 
  (a_1 + a_4)^2 + a_1^2 + a_4^2 + (a_1 + a_3)^2 + a_2^2 + a_3^2
     \,.
     \cr
}
\tag 7.4
$$
We must find those integers $a_1$ and $a_4$ such that
$(a_1 + a_4)^2 + a_1^2 + a_4^2 \leq 2$ and similarly
for $a_2$ and $a_3$.
This yields the list
$$
(a_1,\, a_4) \ \in \
  \left\{
    (1,\, -1),\ (0,\, -1),\  (0,\, 1),\  (1,\, 0),\  (0,\, 0),\
    (-1,\, 0),\ (-1,\, 1)
  \right\}
\tag 7.5
$$
and the same for the pair $(a_2,\, a_3)$.
Putting together those four tuples for which the expression
in (7.4) adds up to two yields the following as the 
coefficients of all possible units.
$$
(a_1,\, a_2,\, a_3,\, a_4 ) \ = \ 
\pm \left\{
\eqalign{
     (1,\, 0,\, 0,\, -1),\ (1,\, 0,\, 0,\, 0),\ (0,\, 0,\, 0,\, 1)
  \cr 
     (0,\, 1,\, -1,\, 0),\ (0,\, 1,\, 0,\, 0),\ (0,\, 0,\, 1,\, 0)
  \cr
}
\right\}
\tag 7.6
$$
which are the units in (7.2).
\qed
\enddemo

\proclaim{Lemma 27}
For every element $\bold a$ in $H_{1,3,3}$ 
there exists two units $\bold u_1, \bold u_2 \in H_{1,3,3}$
such that $\bold u_1 \, 2 \,\bold a \, \bold u_2$ is an
element of 
$\Bbb Z [1,\, \qi,\, \sqrt 3 \,\qj,\, \sqrt 3 \, \qk]$.
\endproclaim

\demo{Proof}
Write 
$\bold a = a_1 \bold x_1 + a_2 \bold x_2 + a_3 \bold x_3 + a_4 \bold x_4$
in the form
$$
\bold a \ = \
  2 \, (c_1 \bold x_1 + c_2 \bold x_2 + c_3 \bold x_3 + c_4 \bold x_4 )
    \,+\,
  d_1 \bold x_1 + d_2 \bold x_2 + d_3 \bold x_3 + d_4 \bold x_4
\tag 7.7
$$
where $d_1, \dots , d_4 \in \{ 0,\, 1 \}$
and $c_1, \dots, c_4$ are rational integers.
Let 
$\bold c = \sum_{t=1}^4 c_t \bold x_t$
and
$\bold d = \sum_{t=1}^4 d_t \bold x_t$.
Computation shows that for every possible value of
$\bold d$ there exist units $\bold u_1$, $\bold u_2$
in $H_{1,3,3}$ such that $\bold u_1 \, 2 \,\bold d \, \bold u_2$ 
is an element of 
$\Bbb Z [1,\, \qi,\, \sqrt 3 \,\qj,\, \sqrt 3 \, \qk]$.
We consider $\bold u_1 \, 2 \, \bold a \, \bold u_2$
and the demonstration concludes as in the previous cases.
\qed
\enddemo

We now proceed to an extension of a trick of Euler
that was used in demonstrating the classical Four
Squares Theorem.

\proclaim{Lemma 28}
Suppose $n$ is an integer for which $2 \, n$ can be
represented by the quadratic form $(1,\, 1,\, 3,\, 3)$.
Then $n$ can also be represented by this form.
\endproclaim

\demo{Proof}
Start with $2 \, n = x^2 + y^2 + 3 z^2 + 3 w^2$ with
$x,\, y,\, z,\, w \in \Bbb Z$.
Consider this modulo $2$.
$$
0 \ \equiv \ 
  x^2 + y^2 + 3 z^2 + 3 w^2
    \ \equiv \
  x + y + z + w \pmod 2
    \,.
\tag 7.8
$$
If $x$ and $y$ have the same parity, then so do $z$ and 
$w$.
We can directly use Euler's trick.
Since $(x+y) / 2$, $(x-y) /2$, $(z + w) / 2$ and $(z - w) / 2$
are integers we compute
$$
\left( {{x+y} \over 2} \right)^2  \, + \,
  \left( {{x-y} \over 2} \right)^2  \, + \, 
  3 \, \left( {{z+w} \over 2} \right)^2 \, + \,
  3 \, \left( {{z-w} \over 2} \right)^2 
    \ = \ 
  {1 \over 2} \, 2 n 
    \ = \
   n
    \,. 
\tag 7.9
$$
The other case is if $x$ and $y$ do not have the same parity.
With no loss of generality, assume $x$ is even and $y$ is odd.
Then 
$1 \equiv z + w \pmod 2$
so $z$ and $w$ must have different parities too.
Again with no loss of generality suppose $z$ is odd and
$w$ is even.
Consider $y^2 + 3 z^2$ where $y$ and $z$ are both odd.
Set
$$
y_1 \ = \ {{y + 3 z} \over 2}, \quad
  z_1 \ = \ {{y - z} \over 2}, \qquad
y_2 \ = \ {{y - 3 z} \over 2}, \quad
  z_2 \ = \ {{y + z} \over 2}
    \,.
\tag 7.10
$$
Then
$$
\eqalign{
  y_1^2 + 3 z_1^2
      \ &= \
    {1 \over 4}\, \left( y^2 + 6 yz + 9 z^2 + 3 (y^2 - 2yz + z^2) \right)
      \cr
      \ &= \ 
    {1 \over 4} \, ( 4y^2 + 12 z^2)
      \ = \ 
    y^2 + 3 z^2
      \,. \cr
}
\tag 7.11
$$
Similarly $y_2^2 + 3 z_2^2 = y^2 + 3 z^2$.

Since $y$ and $z$ are odd they must be congruent to $1$ or $3$ 
modulo $4$.
If $y \equiv z \pmod 4$ then $(y -z) / 2$ is even and so is
$(y + 3z) / 2$.
Hence $y_1$ and $z_1$ are even, but
$$
x^2 + y^2 + 3 z^2 + 3 w^2
    \ = \ 
  x^2 + y_1^2 + 3 z_1^2 + 3 w^2
    \ = \ 
  2n
\tag 7.12
$$
and we are in the case where $x$ and $y_1$ have the same parity.
Thus $n$ is represented by $(1,\, 1,\, 3,\, 3)$ using the 
arguments of the same parity case above.

If $y \not\equiv z \pmod 4$ then $y + z \equiv 0 \pmod 4$
so $(y+z) / 2$ and $(y - 3z) / 2$ are even.
Thus $y_2$ and $z_2$ are even and
$$
x^2 + y^2 + 3 z^2 + 3 w^2
    \ = \ 
  x^2 + y_2^2 + 3 z_2^2 + 3 w^2
    \ = \ 
  2n
  \,.
\tag 7.13
$$
We are in the same parity case and a representation of $n$ by
$(1,\, 1,\, 3,\, 3)$ is derived as above.
\qed
\enddemo

\proclaim{Lemma 29}
$H_{1,3,3}$ is norm Euclidean.
\endproclaim

\demo{Proof}
Fix any quaternion $\bold q \in \Bbb H$.
We first find an integer $a_3$ such that 
$\bold q - a_3 \bold x_3$ has $\qj$ component less than or equal
to $\sqrt 3 / 4$ in absolute value.
We must now consider the $\qi$ and $\qk$ components
simultaneously.
Let $\bold x_2^{\hbox{\sevenrm P}}$ and $\bold x_4^{\hbox{\sevenrm P}}$ 
be the projections of
$\bold x_2$ and $\bold x_4$ onto the $\qi$--$\qk$ plane,
$$
\bold x_2^{\hbox{\sevenrm P}}
  \ = \ 
       {1 \over 2} \, \qi - {{\sqrt 3} \over 2} \, \qk,
    \qquad
  \bold x_4^{\hbox{\sevenrm P}}  
  \ = \ 
      {3 \over 4} \, \qi + {{\sqrt 3} \over 4} \, \qk
    \,.
\tag 7.14
$$
We wish to approximate an arbitrary point in the $\qi$--$\qk$
plane, $a\qi + b \qk$, by an integral linear combination of
the above two projections.
Set $m$ equal to the nearest integer to 
$(a- \sqrt 3 \, b) / 2$ and $n$ equal to the nearest integer to
$(3a + \sqrt 3 \, b) /3$.
Thus
$$
\left| m - {{a - \sqrt 3 \, b} \over 2} \right| 
  \ \leq \
    {1 \over 2},
  \qquad
    \left| n - {{3a + \sqrt 3 \, b} \over 3} \right| 
  \ \leq \
    {1 \over 2}
  \,. 
\tag 7.15
$$
Set $-c$ equal to the expression inside the absolute value
signs on the left, while $-d$ is set equal to the corresponding
term on the right.
Note
$$
a \qi + b \qk - m \bold x_2^{\hbox{\sevenrm P}} 
    - n \bold x_4^{\hbox{\sevenrm P}}
  \ = \ 
    \left(a - {1 \over 2} m - {3 \over 4} n \right) \, \qi
       + \left( b + {{\sqrt 3} \over 2} m - {\sqrt 3 \over 4} n
         \right) \, \qk
\tag 7.16
$$
which implies
$$
|| 
a \qi + b \qk - m \bold x_2^{\hbox{\sevenrm P}} 
    - n \bold x_4^{\hbox{\sevenrm P}}  ||^2
  \ = \ 
    \left( a - {1 \over 2} m - {3 \over 4} n \right)^2
       + \left( b + {{\sqrt 3} \over 2} m - {\sqrt 3 \over 4} n
         \right)^2
\tag 7.17
$$
as $\qi$ and $\qk$ are orthogonal.
We have the further equations
$$
\eqalign{
a - {1 \over 2} m - {3 \over 4} n
  \ &= \ 
    \left( {{a - \sqrt 3 \, b} \over 2} - m \right) \cdot {1 \over 2}
    \ + \
    \left( {{3a + \sqrt 3 \, b} \over 3} -n \right) \cdot {3 \over 4}
  \cr
b + {{\sqrt 3} \over 2} m - {{\sqrt 3} \over 4} n 
  \ &= \ 
    \left( {{a - \sqrt 3 \, b} \over 2} - m \right) \cdot
        \left(  -{{\sqrt 3} \over 2} \right)
    \ + \
    \left( {{3a + \sqrt 3 \, b} \over 3} -n \right) \cdot {{\sqrt 3} \over 4}
  \cr
}
\tag 7.18
$$
Then    
$$
\eqalign{
|| 
a \qi + b \qk - m \bold x_2^{\hbox{\sevenrm P}} 
    - n \bold x_4^{\hbox{\sevenrm P}}  ||^2
  \ &= \ 
    \left( {1 \over 2} \, c + {3 \over 4} \, d \right)^2
       \ + \
    \left( -{{\sqrt 3} \over 2} \, c + {{\sqrt 3} \over 4} \, d \right)^2
  \cr
  \ &= \ 
    c^2 + {3 \over 4} d^2
  \ \leq \
    {7 \over {16}}
}
\tag 7.19
$$
since $c$ and $d$ are bounded by $1/2$ in absolute value.

Hence an integer linear combination of $\bold x_2$ and $\bold x_4$
can be chosen such that the sums of squares of the $\qi$ and $\qk$
components of 
$\bold q - a_3 \bold x_3 - m \bold x_2 - n \bold x_4$
is bounded by $7 / 16$.
The $\qj$ component is unchanged by this transformation of $\bold q$.
We need only choose integral $a_1$ for which 
$\bold q - a_3 \bold x_3 - m \bold x_2 - n \bold x_4 - a_1 \bold x_1$
has real coefficient less than or equal to $1/2$.
Then
$$
||
  \bold q - a_3 \bold x_3 - m \bold x_2 - n \bold x_4 - a_1 \bold x_1
    ||^2 
  \ \leq \
    \left( {1 \over 2} \right)^2
    + \left( {{\sqrt 3} \over 4} \right)^2
    + {7 \over {16}}
  \ = \ 
    {7 \over 8}
  \,.
\tag 7.20
$$
Hence $H_{1,3,3}$ is norm Euclidean.
\qed
\enddemo

\proclaim{Lemma 30}
For any rational prime $p$ there exist rational integers
$a$ and $b$ such that 
$a^2 + b^2 + 3 \equiv 0 \pmod p$.
\endproclaim

\demo{Proof}
The demonstration is analogous to similar Lemmas above.
\qed
\enddemo

\proclaim{Lemma 31}
$H_{1,3,3}$ contains the elements $1$, $2 \, \qi$,
$\sqrt 3 \, \qj$, $\qi + \sqrt 3 \, \qk$.
\endproclaim

\demo{Proof}
$\bold x_1 = 1 \in H_{1,3,3}$.
Also $\bold x_2 + 2 \, \bold x_4 = 1 + 2 \, \qi$ hence
$2 \, \qi \in H_{1,3,3}$.
Then 
$2 \bold x_3 - \bold x_2 = \sqrt 3 \, \qj \in H_{1,3,3}$.
Finally 
$2 \, \bold x_4 - \bold x_2 = 1 + \qi + \sqrt 3 \, \qk$
so $\qi + \sqrt 3 \, \qk \in H_{1,3,3}$.
\qed
\enddemo

\proclaim{Theorem 32}
Every rational prime can be represented by the 
form $(1,\, 1,\, 3,\, 3)$.
\endproclaim

\demo{Proof}
The proof has only a few differences from previous
demonstrations of analogous Theorems.
By similar reasoning to the above mentioned Theorems
we find that there exist integers $a$ and $b$ such that
$a^2 + b^2 +3 = pr$ for integer $r$, $0 < r < p$.

To find an element $\bold u \in H_{1,3,3}$ with norm
equal to $a^2 + b^2 + 3$ we consider two cases.
If $b$ is even we may use
$\bold u = a + (b / 2) \cdot 2 \, \qi + \sqrt 3 \, \qj$.
This is in $H_{1,3,3}$ by the previous Lemma.
If $b$ is odd then write $b = 2 s + 1$ and let
$$
\bold u \ = \ 
  a + s \cdot 2 \, \qi + (\qi + \sqrt 3 \, \qk)
  \ \in \ H_{1,3,3}
  \,.
\tag 7.21
$$
Then $\bold u = a + b \qi + \sqrt 3 \, \qk$ so the norm
of $\bold u$ is again $a^2 + b^2 + 3$.

The remainder of the proof is similar with the note that
it easy to see $H_{1,3,3}$ is closed under conjugation.
We end up with an element $\bold d \in H_{1,3,3}$ such
that $N(\bold d) = p$.
Since our results only guarantee that a two sided
associate of $2 \, \bold d$ is in the ring
$\Bbb Z [1,\, \qi,\, \sqrt 3 \, \qk,\, \sqrt 3 \, \qk]$
we find integers $a_1,\, a_2,\, a_3,\, a_4$ such that
$$
4 \, p \ = \ 
  4 \, N(\bold d) \ = \ 
    N(2 \, \bold d) \ = \ 
      a_1^2 + a_2^2 + 3 \, a_3^2 + 3 \, a_4^2
  \,.
\tag 7.22
$$
By Lemma 28, the form $(1,\, 1,\, 3,\, 3)$
represents $2p$, and hence $p$.
\enddemo

\proclaim{Corollary 33}
The quadratic form $(1,\, 1,\, 3,\, 3)$ is universal.
\endproclaim

\demo{Proof}
This form represents one, every prime, and is norm
quaternionic.
\enddemo

\medskip
%------------------------------------------------
\head 8. The form $(1,\, 2,\, 5,\, 10)$\endhead

The attempt to demonstrate the universality of this
quadratic form via quaternions runs into difficulties.
When looking for a basis of units, we previously used
quaternions with real part $0$, ${1 \over 2}$ and $1$.
Experience shows other real parts lead to dramatically
increasing denominators upon repeated multiplication.
At least one of the unit generators in the above cases
had real part zero.
It is not difficult to show that the unit quaternions
analogous to those considered for the previous quadratic
forms cannot have real part zero.

\proclaim{Theorem 34}
There is no unit quaternion of the form
$$
{1 \over e} \,\, (b\, \sqrt 2 \,\qi + c \,\sqrt 5 \,\qj 
   + d \,\sqrt{10} \,\qk )
\tag 8.1
$$
with $b,\, c,\, d,\, e \in \Bbb Z$.
\endproclaim

\demo{Proof}
For the norm to be one we must have
$e^2 = 2 b^2 + 5 c^2 + 10 d^2$.
With no loss of generality we may choose $e,\, b,\, c,\, d$ to be
positive.
Consider the solution with the smallest positive value of $e$.
Then $e^2 \equiv 2 b^2 \pmod 5$.
If $b \not\equiv 0 \pmod 5$ then $2$ is a quadratic residue modulo
$5$.
This is a contradiction, so $5$ must divide $b$.
Consequently $5$ also divides $e$.

Let $e = 5 e_1$ and $b = 5 b_1$.
Thus $e_1$ and $b_1$ are elements of $\Bbb Z$ and we have
$$
\eqalign{
  25 \, e_1^2 \ &= \ 2 \cdot 25 \, b_1^2 + 5 \, c^2 + 10 \, d^2
    \cr
   5 \, e_1^2 \ &= \ 10 \, b_1^2 + c^2 + 2 \, d^2
    \cr
}
\tag 8.2
$$
Modulo $5$ this becomes $0 \equiv c^2 + 2 d^2 \pmod 5$.
If $d$ is not divisible by $5$ then $-2$ is a quadratic
residue of $5$, a contradiction.
Hence $d$ and thus $c$ are divisible by $5$.

Let $c = 5 c_1$ and $d = 5 d_1$.
Then $c_1$ and $d_1$ are in $\Bbb Z$ and
$$
\eqalign{
   5 \, e_1^2 \ &= \ 10 \, b_1^2 + 25 \, c_1^2 + 2 \cdot 25 \, d_1^2
    \cr
        e_1^2 \ &= \ 2 \, b_1^2 + 5 \, c^2 + 10 \, d^2
    \cr
}
\tag 8.3
$$
which contradicts the minimality of the choice of $e$.
\qed
\enddemo

\medskip
%------------------------------------------------
\head 9. The Computation\endhead
PUNIMAX, a computer algebra system related to MAXIMA, was
used on a LINUX partition of two personal computers.
One computer was a Pentium 133 with 32 megabytes of RAM
and the other was a Pentium 300 with 64 megabytes of RAM.
Both machines ran Linux 2.0.35.
The computation of the two sided products of units for
the quadratic form $(1,\, 2,\, 3,\, 6)$ required a
total of an hour of machine time on the faster computer.
A UNIX port of SNOBOL was used to change MAXIMA output
into a \TeX\ friendly form.
In addition Tables II and IV were checked by manual calculations.
On this basis, Table III and the cases listed in Table V
were also verified manually.
Table VI was verified by an elementary {\sl lex\/} and {\sl yacc\/}
based quaternion calculator written by the author.

\medskip
%------------------------------------------------
\head 10. Acknowledgments\endhead
The author would like to thank A.~Jarvis and N.~Dummigan
of Sheffield University
for assistance in tracking down one of the references and 
explaining some current results.
The author would like to thank B.~Haible, the maintainer
of PUNIMAX, for permitting its free use for academic purposes 
\localcite{Haible}{5}{1997}.
The author would like to thank Harvey Cohn for many
encouraging email communications.
The author expresses his appreciation to the referee for
very useful comments and suggestions.

\refstyle{C}

\enddocument